\documentclass[10pt]{article}
\usepackage{ifpdf}
\usepackage{setspace}
\usepackage{graphicx,amssymb,lineno}
\usepackage{amsmath}
\ifpdf
\usepackage[%
   pdfauthor={Simon Pepping},%
  pdfstartview=FitH,%
  bookmarks=true,%
  bookmarksopen=true,%
  breaklinks=true,%
  colorlinks=true,%
  linkcolor=blue,anchorcolor=blue,%
  citecolor=blue,filecolor=blue,%
  menucolor=blue,pagecolor=blue,%
  urlcolor=blue]{hyperref}
\else
\usepackage[%
  breaklinks=true,%
  colorlinks=true,%
  linkcolor=blue,anchorcolor=blue,%
  citecolor=blue,filecolor=blue,%
  menucolor=blue,pagecolor=blue,%
  urlcolor=blue]{hyperref}
\fi

\makeatletter
\def\elsartstyle{%
    \def\normalsize{\@setfontsize\normalsize\@xiipt{14.5}}
    \def\small{\@setfontsize\small\@xipt{13.6}}
    \let\footnotesize=\small
    \def\large{\@setfontsize\large\@xivpt{18}}
    \def\Large{\@setfontsize\Large\@xviipt{22}}
    \skip\@mpfootins = 18\p@ \@plus 2\p@
    \normalsize
} \@ifundefined{square}{}{} \makeatother

\def\be{\begin{equation}}
\def\ee{\end{equation}}

\begin{document}
\title{A Partial Hamiltonian Approach for Current Value Hamiltonian Systems}
\begin{center}{\large\bf R. Naz$^a$, F. M. Mahomed$^b$ and Azam Chaudhry$^c$}

{$^a$ Centre for Mathematics and Statistical Sciences,
  Lahore School of Economics, Lahore, 53200, Pakistan\\
  $^b$ Differential Equations,
 Continuum Mechanics and Applications, School of Computational
and Applied Mathematics, University of the Witwatersrand, Wits 2050,
South Africa\\
 $^c$ Department of Economics,
 Lahore School of Economics, Lahore, 53200, Pakistan}

\end{center}
{\bf Abstract}\\
We develop a partial Hamiltonian framework to obtain reductions and
closed-form solutions via first integrals of current value
Hamiltonian systems of ordinary differential equations (ODEs).
 The approach is algorithmic and applies to many state and costate
variables of the current value Hamiltonian. However, we apply the
method to models with one control, one state and one costate
variable to illustrate its effectiveness. The current value
Hamiltonian systems arise in economic growth theory and other
economic models.  We explain our approach with the help of a simple
illustrative example and then apply it to two widely used economic
growth models: the Ramsey model with a constant relative risk
aversion (CRRA) utility function and Cobb Douglas technology and a
one-sector AK model of endogenous growth are considered.  We show
that our newly developed systematic approach can be used to deduce
results given in the literature and also to find new solutions.

{\bf keyword:} Current value Hamiltonian, partial Hamiltonian
approach, economic growth models

\section{Introduction}
There has been extensive use of dynamic optimization in economic
modeling and  many of these models use the current value Hamiltonian
whenever the integrand function contains a discount factor. These
models range from those used for neoclassical economic growth
(\cite{ramsey}, \cite{lucas}) to optimal firm-level investment
\cite{eisner} and human capital and earnings \cite{ben}. Pontrygin's
maximum principle provides a set of necessary conditions for the
solution of the continuous time optimal control problem involving a
current value Hamiltonian and a dynamical system of ODEs is obtained
for control, state and costate variables. Beginning with \cite{cass}
there have been various approaches, both qualitative and
quantitative (see \cite{rit} for a good account of these), to deal
with dynamic economic models arising from current value Hamiltonian
system and most of these models were solved using
numerical approaches (like \cite{mul}) or linear approximations
around steady states (\cite{baro}). The critical problem is that for
the underlying nonlinear dynamical system in economics there is a
lack of a general analytical solution procedure not only for higher
order systems but even for systems with one state and costate
variable.

It is true to say that nonlinear dynamical systems evade closed-form
solutions in general. However, the lack of a general procedure
inhibits the search for reductions and solutions of such type of
nonlinear equations even when solutions do exist. Having said that,
there are some well-known closed-form solutions that appear in the
literature (see, e.g. \cite{s1, s2, s3, s4, s5, s6}). These
solutions have been obtained by seemingly disparate approaches.
Independent of the knowledge of explicit solutions, dynamic local
stability of certain systems (see \cite{st1,st2,st3}) have been
characterized by qualitative or numerical approaches.

Several important contributions have been made in the analysis of
nonlinear dynamical systems of economic models. Here we focus on a
new approach which yields first reductions and closed-form solutions
for such systems of ODEs. We develop a Hamiltonian framework for
several control, state and costate variables. Therefore, the method
we develop is applicable to an arbitrary system of ODEs. However, we
apply it to a system of two ODEs in order to show its effectiveness.
In the case of higher order systems of ODEs, the approach may
require the use of algebraic computing.

The layout of the paper is as follows. The partial Hamiltonian
approach is developed in Section 2.  In Section 3 we provide a
simple illustrative example to show how our approach works. The
Ramsey model with a constant relative risk aversion (CRRA) utility
function with Cobb Douglas technology and the one-sector AK model of
endogenous growth are studied in Section 4 and known solutions are
deduced via our partial Hamiltonian approach. Conclusions are
finally presented in Section 5.

\section{A Hamiltonian version of the Noether-type theorem}

Herein we develop a partial Hamiltonian approach for current value
Hamiltonians which do not satisfy the canonical Hamilton equations.
This is done for several control, state and costate variables.

Let $t$ be the independent variable and $(q,p)=(q^1,...,q^n,
p_1,...,p_n)$ the phase space coordinates. The derivatives of $q^i$,
$p_i$ with respect to $t$ are \be \dot p_i=D(p_i),\; \dot
q_i=D(q_i),\;i=1,2,\cdots,n, \label{(1)} \ee where \be
D=\frac{\partial }{\partial t}+\dot q_i\frac{\partial }{\partial
q_i}+\dot p_i\frac{\partial }{\partial p_i}+\cdots\label{(2)} \ee is
the total derivative operator with respect to $t$. The summation
convention is utilized for repeated indices. The variables $t, q, p$
are independent and connected only by the differential relations
{(\ref {(1)})}.

There are some well-known operators which are defined in the space
of the variables $(t,q,p)$ and its prolongations. We introduce them.

In addition to the Euler operator \be \frac{\delta}{\delta
q^i}=\frac{\partial}{\partial q^i}-D\frac{\partial}{\partial \dot
q^i}, i=1,2, \cdots, n, \label{(3)}\ee one also has the variational
operator
 \be \frac{\delta}{\delta
p_i}=\frac{\partial}{\partial p_i}-D\frac{\partial}{\partial \dot
p_i}, i=1,2, \cdots, n. \label{(4)}\ee The action of the operators
{(\ref {(3)})} and {(\ref {(4)})} on
 \be L(t,q,\dot q)=p_i \dot
q^i-H(t,q,p)\label{(5)} \ee equated to zero yields the canonical
Hamilton equations
 \begin{eqnarray} \dot q^i=\frac{\partial H}{\partial p_i},
\nonumber\\
[-1.5ex]\label{(6)} \\[-1.5ex]
 \dot p_i=-\frac{\partial H}{\partial q^i}, i=1,
\ldots, n.  \nonumber
\end{eqnarray}
 That is ${\delta L}/{\delta q^i}=0$ and ${\delta L}/{\delta
 p_i}=0$ results in {(\ref {(6)})}. Equation {(\ref {(5)})} is the
 well-known Legendre transformation which relates the Hamiltonian
 and Lagrangian, where $p_i=\partial L/\partial \dot q^i$ and
$\dot q^i=\partial H/\partial \dot p_i$.

 Generators of point symmetries in the space $(t,q,p)$ are operators of the form
 \be X=\xi(t,q,p)\frac{\partial }{\partial t}+ \eta ^i (t, q,p)\frac{\partial }{\partial q^i}
 + \zeta _i (t,q,p) \frac{\partial }{\partial p_i}. \label{(7)}\ee
 The operator in (\ref{(7)}) is a generator of a point symmetry of the canonical Hamiltonian system
 (\ref{(6)}) if (\cite{olv})
  \begin{eqnarray}
  \dot \eta^i-\dot q^i\dot\xi-X(\frac{\partial H}{\partial p_i})=0,\nonumber\\
  \dot\zeta_i-\dot p_i\dot\xi+X(\frac{\partial H}{\partial q^i})=0,\; i=1,\ldots,n \label{(7a)}
  \end{eqnarray}
  on the system (\ref{(6)}).

 Hamiltonian symmetries in evolutionary or canonical form have been
 considered (\cite{olv}).
 Furthermore, symmetry properties of the
 Hamiltonian action have been investigated in the space $(t,q,p)$ by
\cite{leach}  and  \cite{roman}. In the latter, the authors
considered
 the general form of the symmetries {(\ref {(7)})} and provided a
 Hamiltonian version of Noether's theorem.

 The following important results which are analogs of Noether symmetries and the Noether
 theorem (see \cite{olv, ibra1, ibra2, roman} for a discussion) were established.

 {\bf Theorem 1} (Hamilton action symmetries):
 A Hamiltonian action
 \be p_idq^i -H dt  \label{(8)}\ee
 is invariant up to gauge $B(t, q,p)$  with respect to a group generated
 by {(\ref {(7)})} if and only if the condition
 \be \zeta_i \frac{\partial H}{\partial p_i} + p_i D(\eta ^i)-X(H)-H D(\xi)-D(B)=0, \label{(9)}\ee
 holds.

 { \bf Theorem 2} (Hamiltonian version of Noether's theorem):
The canonical Hamilton system {(\ref {(6)})} which is invariant has
the first integral
 \be I= p_i \eta ^i -\xi H-B \label{(10)}\ee
 for some gauge function $B=B(t,q,p)$ if and only if the Hamiltonian
 action is invariant up to divergence with respect to the operator $X$ given in {(\ref {(7)})}
on the solutions to equations (\ref{(6)}).

We now focus our attention on systems of equations which are not in
the canonical form (\ref{(6)}). Therefore the Theorems 1 and 2 do
not apply for these systems. We need an extension of the existing
results which we carry out below.

 Since the current value Hamiltonian (see, e.g. \cite{chiang}) satisfies
\begin{eqnarray} \dot q^i=\frac{\partial H}{\partial p_i},
\nonumber\\
[-1.5ex]\label{(11)} \\[-1.5ex]
 \dot p^i=-\frac{\partial H}{\partial q^i}+\Gamma_i,\; i=1, 2,
\cdots, n,  \nonumber
\end{eqnarray}
where $\Gamma_i$ is a nonzero function, we seek an extension of the
results relating to the canonical Hamiltonian system to the system
 {(\ref{(11)})} so that we can obtain first integrals of system {(\ref
{(11)})} in an algorithmic manner. We also refer to an $H$ that
satisfies (\ref{(11)}) as a partial Hamiltonian.

It is opportune to remark that $X$ as in (\ref{(7)}) is a generator
of point symmetry of the current value Hamiltonian system
(\ref{(11)}) if
\begin{eqnarray}
  \dot \eta^i-\dot q^i\dot\xi-X(\frac{\partial H}{\partial p_i})=0,\nonumber\\
  \dot\zeta_i-\dot p_i\dot\xi+X(\frac{\partial H}{\partial q^i}-\Gamma_i)=0,\; i=1,\ldots,n \label{(11a)}
  \end{eqnarray}
  on the system (\ref{(11)}). Note that (\ref{(11a)}) is evidently different from (\ref{(7a)}) due to the
  nonzero term $\Gamma_i$.

We introduce the definition of what we call the partial Hamiltonian
operator below. This is motivated by the analogous definition of the
partial Noether operator given in \cite{{imran},{kar3}}.

{\bf Definition 1}: An operator $X$ of the form {(\ref{(7)})} is a
partial Hamiltonian operator corresponding to a current value
Hamiltonian as in {(\ref{(11)})}, if there exists a function
$B(t,q,p)$ such that

 \be \zeta_i \frac{\partial H}{\partial p_i}+p_i D(\eta^i)-X(H)
-H D(\xi)=D(B) +(\eta^i -\xi \frac{\partial H}{\partial
p_i})(-\Gamma_i) \label{(12a)}\ee holds.

Note that if $H$ is a present value Hamiltonian, then equation
{(\ref{(12a)})} becomes the usual determining equation for
symmetries of the Hamiltonian action since $\Gamma_i=0$ in this
case.

Also one can immediately see from {(\ref{(12a)})} that if $X$ and
$Y$ are partial Hamiltonian operators, then so is a linear
combination of these Hamiltonian operators.

We now have the following important theorem on how one constructs
first integrals for the system  {(\ref{(11)})}. That is we present
the partial Hamiltonian approach for current value Hamiltonians.
This is achieved for several control, state and costate variables.

{\bf Theorem 3} (partial Hamiltonian version pf the partial Noether
theorem): An operator $X$ of the form  {(\ref{(7)})} is a partial
Hamiltonian operator of the current value Hamiltonian $H$
corresponding to system  {(\ref{(11)})} if and only if
{(\ref{(10)})} is its first integral.

{\bf Proof}: The result follows by straightforward differentiation
of the first integral formula  {(\ref{(10)})} on the solutions of
system {(\ref{(11)})}. However, one has to remember that the terms
involving $\dot p_i$ need to be replaced by the right hand side of
the second equation of system {(\ref{(11)})} which has a non-zero
function $\Gamma_i$. Another way to show this to be the case is to
utilize \be D(H)|_{\dot q^i=\frac{\partial H}{\partial p_i}, \dot
p^i=-\frac{\partial H}{\partial q^i}
+\Gamma_i}=H_t+\Gamma_i\frac{\partial H}{\partial p_i}\label{a} \ee
as well as the identity
\begin{eqnarray}
\zeta_i \dot q^i+p_i D(\eta^i)-X(H) -H D(\xi)-D(B) -(\eta^i -\xi
\frac{\partial H}{\partial
p_i})(-\Gamma_i)\nonumber\\
=\xi(D(H)-H_t-\Gamma_i\frac{\partial H}{\partial p_i})-
\eta^i(\dot p^i+\frac{\partial H}{\partial q^i}-\Gamma_i)+\zeta_i(\dot q^i-\frac{\partial H}{\partial p_i})\nonumber\\
+D(p_i \eta ^i -\xi H-B) \label{b}
\end{eqnarray}
which holds for any smooth functions $H(t,q,p)$ and suitable
functions $B(t,q,p)$ and $\Gamma_i$. This identity follows from
direct computations.

{\bf Remark}. An approach in proving Theorem 3 is by invoking the
Legendre transformation (\ref{(5)}) on the partial Noether operators
and partial Noether theorem given in \cite{imran}.

\section{A Simple illustrative example}
Consider the following mathematical example:

 Maximize \be \int _0^{\infty} [\alpha q-\beta q^2- \alpha u^2-
\gamma u] e^{-rt}dt\label{(12)}\ee subject to \be \dot q=u,
\label{(13)}\ee where $\alpha, \beta, \gamma$ are all positive, $r$
is a discount factor, $q(t)$ is the state variable and $u(t)$ is the
control variable.\\\\
{\it \large Hamiltonian function and maximum principle}:\\
The current value Hamiltonian function is defined as \be
H(t,q,p,u)=\alpha q-\beta q^2- \alpha u^2- \gamma u+pu
\label{(14)}\ee where $p(t)$ is called the costate variable. The
necessary first order conditions for optimal control are
\cite{chiang}: \be \frac{\partial H}{\partial u}=0 \label{(15)} \ee
\be \dot q=\frac{\partial H}{\partial p} \label{(16)} \ee
 \be \dot
p=-\frac{\partial H}{\partial q}+r p \label{(17)} \ee Equation
{(\ref{(15)})}-{(\ref{(17)})} with $H$ given by {(\ref {(14)})}
yields
 \be
p=2\alpha u+\gamma \label{(18)} \ee \be   \dot q= u \label{(20)} \ee
\be \dot p= 2\beta q-\alpha+pr \label{(19)} \ee

 Equations {(\ref{(18)})}-{(\ref{(19)})} need to be solved for
 $p(t), q(t), u(t)$. Of course, the direct way to solve this problem is to eliminate
 $p$, $u$ by utilizing {(\ref{(18)})}-{(\ref{(19)})} in order to obtain a scalar linear second order
 ordinary differential equation in $q$, which is amenable to straightforward integration.
 We  explain here how we can find the solution
 by using the partial Hamiltonian approach introduced above.\\\\
{\it \large Determination of Partial Hamiltonian operators}:\\
The partial Hamiltonian operator determining equation  is given in
{(\ref{(12a)})} Expansion of equation {(\ref {(12a)})} yields
 \begin{eqnarray}p(\eta_t+\dot q \eta _q)-\eta (\alpha-2\beta q)-(\alpha q-\beta q^2- \alpha u^2- \gamma u+pu)(\xi_t+\dot q \xi_q)
 \nonumber\\
 [-1.5ex]\label{(22)} \\[-1.5ex]
 =B_t+\dot q B_q +(\eta-\xi u)(-r p),\nonumber
\end{eqnarray}
in which  we assume that $\xi=\xi (t,q) $, $\eta=\eta (t,q) $,
$B=B(t,q)$.

Note that one can also assume these functions to be dependent on
$p$. We have chosen $(t,q)$ dependence to simplify the calculations
here and in the models considered in Section 4 with the purpose of
deriving solutions. This assumption leads to at least one partial
Hamiltonian operator. More general assumptions are used in the event
that one does not obtain an operator.

 With the help of {(\ref {(18)})}-{(\ref
 {(20)})},  Equation {(\ref {(22)})} can be written as
\begin{eqnarray} (2\alpha u+\gamma)(\eta_t+u \eta_q)
 -(\xi_t+u\xi_q)(\alpha q-\beta q^2-\alpha u^2-\gamma u+2 \alpha u^2+\gamma u)
\nonumber\\
[-1.5ex]\label{(23)} \\[-1.5ex]
-\eta(\alpha-2\beta q) =B_t+uB_q+(\eta -\xi u)(-2r \alpha u-\gamma
r). \nonumber
\end{eqnarray}
Separating equation {(\ref {(23)})} with respect to powers of $u$ as
$\xi, \eta, B$ do not contain $u$, we have
\begin{eqnarray} u^3: -\alpha \xi_q=0,\label{(24)} \\
u^2: 2\alpha \eta_q-\alpha \xi_t=2\alpha r \xi,
\label{(25)}\\
u: 2\alpha \eta_t+\gamma \eta_q=B_q-2\alpha r \eta +\gamma r \xi,
\label{(26)}\\
u^0: \gamma \eta_t-\eta(\alpha-2\beta q)-\xi_t (\alpha q- \beta
q^2)=B_t-r \gamma \eta. \label{(27)}
\end{eqnarray}
System {(\ref {(24)})}-{(\ref {(26)})} yields
\begin{eqnarray}\xi=a(t), \; \eta=(\frac{1}{2}
\dot a+r a)q+b(t),
\nonumber\\
[-1.5ex]\label{(28)} \\[-1.5ex]
B=\alpha (\frac{1}{2} \ddot a +r\dot a)q^2+\alpha r (\frac{1}{2}
\dot a +r a)q^2+2\alpha\dot b q+2 \alpha r b q+\frac{1}{2} \gamma
\dot a q+d(t). \nonumber
\end{eqnarray}
Substituting $\xi, \eta, B$ from {(\ref {(28)})} in {(\ref {(27)})}
and then separating w.r.t powers of $q$ we have

\begin{eqnarray}q^2:\frac{1}{2}\alpha \dddot a+\frac{3}{2}\alpha r
\ddot a+(\alpha r -2\beta )\dot a-2\beta r a=0, \qquad
\label{(29)}\\
q: \frac{3}{2}(r\gamma-\alpha)\dot a+r
(r\gamma-\alpha)a+2b\beta=2\alpha \ddot b+2\alpha r \dot b, \label{(30)}\\
q^0: \gamma \dot b-\alpha b+\gamma r b=\dot d. \qquad
\qquad\qquad\qquad\qquad\qquad\label{(31)}
\end{eqnarray}
The solution of equations {(\ref {(29)})}-{(\ref {(31)})} for $a,
b,d$ with general $\alpha, \beta, \gamma, r$ is purely formal and
depend on the roots of the characteristic equation. Clearly there
are three lengthy solutions for $a$ and two for $b$. To be
transparent, we have selected values. Therefore we seek a solution
of equations {(\ref {(29)})}-{(\ref {(31)})} for $a, b,d$ with
specific values $\alpha=\gamma=r=1$ and $\beta=2$ and we arrive at
\begin{eqnarray}
a(t)=c_1e^{-t}+c_2e^{2t}+c_3e^{-4t}, \nonumber\\
b(t)=c_4e^t+c_5e^{-2t},\nonumber\\
[-1.5ex]\label{(32)} \\[-1.5ex]
d(t)=c_4e^t+c_5e^{-2t}+c_6, \nonumber
\end{eqnarray}
where $c_1,\cdots, c_6$ are arbitrary constants. Finally, we obtain
the following $\xi, \eta, B$ after substituting $a,b,d$ from {(\ref
{(32)})}  into {(\ref {(28)})}
\begin{eqnarray}
\xi=c_1e^{-t}+c_2e^{2t}+c_3e^{-4t}, \nonumber\\
\eta=(\frac{1}{2}c_1e^{-t}+2c_2e^{2t}-c_3e^{-4t})q+c_4e^t+c_5e^{-2t},\nonumber\\
[-1.5ex]\label{(33)} \\[-1.5ex]
B(t)=(6c_2e^{2t}+3c_3e^{-4t})q^2+(-\frac{1}{2}c_1e^{-t}+c_2e^{2t}-2c_3e^{-4t}\nonumber\\
+4c_4e^t-2c_5e^{-2t})q+c_4e^t+c_5e^{-2t}+c_6. \nonumber
\end{eqnarray}
The generators $X_i$ form a vector space. By choosing one of the
constants as one and the rest as zero in turn we have the following
five operators and gauge terms:
\begin{eqnarray}
 X_1=e^{-t}\frac{\partial}{\partial t}+\frac{1}{2}qe^{-t}\frac{\partial}{\partial q},\; B_1=-\frac{1}{2}e^{-t}q\nonumber\\
 X_2=e^{2t}\frac{\partial}{\partial t}+2qe^{2t}\frac{\partial}{\partial q},\; B_2=6q^2e^{2t}+qe^{2t}\nonumber\\
 X_3=e^{-4t}\frac{\partial}{\partial t}-qe^{-4t}\frac{\partial}{\partial q},\; B_3=3q^2 e^{-4t}-2qe^{-4t}\nonumber\\
[-1.5ex]\label{(34)} \\[-1.5ex]
X_4=e^t\frac{\partial}{\partial q},B_4=4qe^t+e^t\nonumber\\
X_5=e^{-2t}\frac{\partial}{\partial q}, B_5=-2qe^{-2t}+e^{-2t}.
\nonumber
\end{eqnarray}
In general the $X_i$'s are not symmetries of the system
\begin{eqnarray}
\dot q=\frac12 p-\frac12,\nonumber\\
\dot p=4q+p-1,\label{d}
\end{eqnarray}
which considered as a non-canonical Hamiltonian system
{(\ref{(18)})}-{(\ref{(19)})} admits the operators (\ref{(34)}). For
example in the case of $X_4$ we have that the first of equations
(\ref{(11a)}) gives $\zeta=2e^t$. However, the second equation of
(\ref{(11a)}) is not
satisfied as easily can be verified.\\\\
{\it \large Construction of first integrals from partial Hamiltonian operators and gauge terms}:\\
Now, first integrals satisfying $DI=0$, on the solutions,
corresponding to operators and gauge terms given in {(\ref {(34)})}
can be computed from {(\ref {(10)})} and the following integrals
result.
\begin{eqnarray} I_1=  [\frac{1}{2}pq   -( q-2q^2-  u^2- u+pu)+\frac{1}{2}]e^{-t},\nonumber\\
I_2=  [2pq  -( q-2q^2-  u^2- u+pu)-6q^2-q]e^{2t}, \nonumber\\
[-1.5ex]\label{(35)} \\[-1.5ex]
I_3=  [-pq   -( q-2q^2-  u^2- u+pu)-3q^2+2q]e^{-4t},\nonumber\\
I_4=  [p -4q-1]e^t,\nonumber\\
I_5=[p+2q-1]e^{-2t} .\nonumber
\end{eqnarray}
There are five first integrals, two of which are functionally
independent.\\\\
{\it \large Optimal solution via first integrals}:\\
Equations {(\ref{(18)})}-{(\ref{(19)})}
 need to be solved for
 $p(t), q(t), u(t)$ with $\alpha=\gamma=r=1, \beta=2$. We demonstrate here how one can find a solution
 by using first integrals. We  derive
the solution associated with the first integral $I_4$. As $DI=0$, on
the solutions, and thus $I=constant$, we have \be [p -4q-1]e^t=A_1,
\label{(36)}\ee where $A_1$ is an arbitrary constant and this gives
 \be p(t) =4q+1+A_1e^{-t}.\label{(37)}\ee
From  {(\ref{(18)})}, $u=\frac{p-1}{2}$ and after using $p$ from
{(\ref{(37)})}, we have
 \be u(t)=\frac{4q+A_1e^{-t}}{2}. \label{(38)}\ee
Thus if $q(t)$ is known we can get the optimal path $p(t)$ and
$u(t)$ from {(\ref{(37)})} and {(\ref{(38)})}.
 Equation {(\ref{(20)})} with $u$ from {(\ref{(38)})} yields
  \be   \dot q= \frac{4q+A_1e^{-t}}{2},\label{(39)}\ee
and this is a first order linear equation in $q(t)$. The solution of
equation {(\ref{(39)})} is
 \be q(t)=\frac{A_1}{2}e^{-t}+A_2e^{2t},\label{(40)} \ee
 where $A_1$ and $A_2$ are arbitrary constants which we can specify
 if we have given initial and terminal conditions.
 One can use any one of the first integrals {(\ref{(35)})} to obtain
 the general solution to this linear problem. Generally a first integral
 provides a reduction of order of the system by one. One can also achieve
 reduction to quadratures of the system by invoking first integrals as we see
 in the examples that follow.

\section{Optimal path of some economic models}
\subsection{Ramsey neoclassical model with CRRA utility function} We consider the following Ramsey
 neoclassical growth model
\cite{baro}, \cite{runge}, where the representative consumer's
utility maximization problem is  defined as
 \be {Max} \quad \int_0^{\infty} e^{-r t}c^{1-\sigma} dt ,\; \sigma \not=0,1\label{(e1)}\ee
subject to the capital accumulation equation and parameter
restriction \be \dot k(t) =k^\beta -\delta k-c, \; k(0)=k_0,\;
0<\beta<1, \label{(e2)}\ee where $c(t)$ is the consumption per
person, $k(t)$ is the capital labor ratio, $\beta, \delta, r$ are
the capital share, depreciation rate, rate of time preferences
respectively. The intertemporal elasticity of substitution is given
by $1/{\sigma}$ and $k_0$ is the initial capital stock.

The current value Hamiltonian function for this model is defined as
\be H(t,c,k,\lambda)=c^{1-\sigma}+\lambda(k^\beta -\delta k-c)
\label{(e3)},\ee where $\lambda(t)$  is the costate variable. The
necessary first order conditions for optimal control are
 \be
\lambda=(1-\sigma)c^{-\sigma} \label{(e4)} \ee
 \be
\dot k=k^\beta -\delta k-c \label{(e5)} \ee
 \be   \dot
\lambda= -\lambda(\beta k^{\beta-1}-\delta)+\lambda r \label{(e6)}
\ee  and the transversality condition is  \be \lim_{t\to\infty}
e^{-rt}\lambda(t)k(t)=0  \label{(trans)}. \ee From {(\ref {(e4)})}
and {(\ref {(e6)})}, the growth rate of consumption is given by \be
\frac{\dot c}{c}= \frac{\beta}{\sigma}
k^{\beta-1}-\frac{1}{\sigma}(\delta+r). \label{(con)} \ee
 We seek a solution  $\lambda(t), k(t), c(t)$ of equations {(\ref{(e4)})}-{(\ref{(e6)})}
 by utilizing the Hamiltonian approach.
The partial Hamiltonian determining equation
 {(\ref{(12a)})}  for the Hamiltonian {(\ref{(e3)})} yields
\begin{eqnarray} \lambda(\eta_t+\dot k \eta _k)-\eta \lambda(\beta k^{\beta-1}-\delta)
 -[c^{1-\sigma}+\lambda(k^\beta -\delta k-c)](\xi_t+\dot k \xi_k)
 \nonumber\\
 [-1.5ex]\label{(e7)} \\[-1.5ex]
 =B_t+\dot k B_k +(\eta-\xi {\partial H\over\partial \lambda})(-r \lambda) ,\nonumber
\end{eqnarray}
 in which we assume that $\xi=\xi (t,k) $,
$\eta=\eta (t,k) $, $B=B(t,k)$. The same reason applies here as for
the illustrative example.
 Equation {(\ref {(e7)})}  with the help of {(\ref {(e4)})}-{(\ref
 {(e6)})} can be written as
\begin{eqnarray}
(1-\sigma)c^{-\sigma}[\eta_t+(k^{\beta}-\delta k-c)\eta_k]-\eta (1-\sigma)c^{-\sigma}(\beta k^{\beta -1}-\delta)\nonumber\\
[-1.5ex]\label{(e8)} \\[-1.5ex]
 -[c^{1-\sigma}+(1-\sigma)c^{-\sigma}(k^\beta -\delta k-c)][\xi_t+(k^\beta -\delta k-c)\xi_k]
 \nonumber\\
 =B_t+(k^\beta -\delta k-c)B_k-r(1-\sigma)c^{-\sigma}[\eta-\xi(k^\beta -\delta
 k-c)].\nonumber
\end{eqnarray}
Separating equation {(\ref {(e8)})} with respect to powers of the
control variable $c$, we have
\begin{eqnarray} c^{2-\sigma}: -\sigma \xi_k=0,\label{(e9)} \\
c^{1-\sigma}: -\eta _k(1-\sigma)-\sigma \xi_t+r(1-\sigma)\xi=0,
\label{(e10)}\\
c^{-\sigma}:\eta_t+(k^{\beta}-\delta k)\eta_k-\eta (\beta
k^{\beta-1}-\delta)\nonumber\\
-(k^{\beta}-\delta k)\xi_t+r \eta -r \xi(k^{\beta}-\delta k)=0 ,
\label{(e11)}\\
c ,\; c^0: B_k=0,\;B_t=0.  \label{(e12)}
\end{eqnarray}
Equations  {(\ref {(e9)})}, {(\ref {(e10)})} and {(\ref {(e12)})}
result in \be \xi=a_1(t),\; \eta=(-\frac{\sigma}{1-\sigma} \dot
a_1+r a_1)k+a_2(t), \; B=0. \label{(e13)} \ee
 Equation {(\ref{(e11)})} with $\xi, \eta, B$ from {(\ref {(e13)})} gives $a_2=0$
 and then reduces to
\begin{eqnarray}
k^{\beta}:  \dot a_1-\frac{\beta r(1-\sigma)}{\beta \sigma -1}a_1=0,\; \beta\sigma\not=1,\label{(e14)}\\
k:-\sigma \ddot a_1+[r(1-2\sigma)+\delta(1-\sigma)]\dot
a_1+r(1-\sigma)(r+\delta)a_1=0.\label{(e15)}
\end{eqnarray}
 Equation {(\ref{(e14)})} is valid if $\sigma \beta\not=1$, for
the case where the capital's share is not equal to the intertemporal
elasticity of substitution. Equations {(\ref{(e14)})} and
{(\ref{(e15)})} yield
 \be a_1(t)=c_1e^{\delta \beta (1-\sigma)t} \label{(e16)}\ee
 with
 \be \sigma=\frac{r+\delta}{\beta \delta} .\label{(e17)}\ee
The restriction on the parameters  {(\ref{(e17)})} is the same as
the one given in \cite{{baro},{runge}} and our approach yields this
during the solution process. Now $\xi, \eta$ and $ B$ are given
by\be \xi= c_1e^{\delta \beta (1-\sigma)t},\; \eta = -c_1\delta
e^{\delta \beta (1-\sigma)t}k, B=0,\label{(ee17)} \ee and the only
partial Hamiltonian operator is \be X=e^{\delta \beta
(1-\sigma)t}\frac{\partial}{\partial t}-\delta e^{\delta \beta
(1-\sigma)t}k\frac{\partial}{\partial k},\; B=0.\label{(e18)} \ee
 The following first integral corresponding to the partial Hamiltonian operator
and gauge terms given in {(\ref {(e18)})} can be computed from
{(\ref {(10)})}: \be I= e^{\delta \beta (1-\sigma)t}[-\sigma
c^{1-\sigma}+(\sigma-1)c^{-\sigma }k^{\beta}].\label{(e19)} \ee We
write {(\ref {(e19)})} as a constant, i.e. \be  -\sigma
c^{1-\sigma}+(\sigma-1)c^{-\sigma }k^{\beta}=A_1e^{\delta \beta
(\sigma-1)t}.\label{(e20)} \ee
 From equation {(\ref {(e20)})}, we have \be
k=[\frac{A_1}{\sigma-1}c^{\sigma}e^{\delta \beta
(\sigma-1)t}+\frac{\sigma}{\sigma-1}c]^{\frac{1}{\beta}}.\label{(e22)}\ee
Our next goal is to get either $c$ or $k$. If $A_1=0$ we arrive at
the well-known solution given in \cite{{baro},{runge}}. Equation
{(\ref {(e22)})} for $A_1=0$ yields \be c(t)=(1-\frac{\beta
\delta}{r+\delta})k^{\beta}\label{(e30)}\ee where
$\frac{\sigma-1}{\sigma}=1-\frac{\beta \delta}{r+\delta}$ by
{(\ref{(e17)})}. Substituting $c$ from equation {(\ref {(e30)})} in
Equation {(\ref {(e5)})} results in
 \be \dot k+\delta k=(\frac{\beta
\delta}{r+\delta})k^{\beta}. \label{(e31)} \ee The solution of
equation {(\ref {(e31)})} subject to the initial condition
$k(0)=k_0$ is given by \be
k(t)=[\frac{\beta}{r+\delta}+(k_0^{1-\beta}-\frac{\beta}{r+\delta})e^{-(1-\beta)\delta
t}]^{\frac{1}{1-\beta}}.\label{(e32)}\ee The solutions {(\ref
{(e30)})} and {(\ref {(e32)})} are the same as the ones derived in
\cite{{baro},{runge}} and satisfy the transversality condition given
by {(\ref {(trans)})}. This guarantees that our approach works. For
$A_1\not=0$, we can get more solutions.
 Now we substitute
{(\ref {(e22)})} into equation {(\ref {(con)})} determining $c$,
viz. \be \frac{d}{dt}(ce^{\beta \delta
t})=\frac{\beta}{\sigma}ce^{\delta \beta
t}[\frac{A_1}{\sigma-1}c^{\sigma}e^{\delta \beta
(\sigma-1)t}+\frac{\sigma}{\sigma-1}c]^{1-\frac{1}{\beta}}.
\label{(e23)}\ee Introducing \be S=ce^{\beta \delta t},\ee equation
{(\ref {(e23)})} directly results in \be
\frac{\beta}{\sigma}e^{\delta(1-\beta)
t}dt=\frac{dS}{S[\frac{A_1}{\sigma-1}S^{\sigma}+\frac{\sigma}{\sigma-1}S]^{1-\frac{1}{\beta}}},\label{(e24)}\ee
which provides the general solution. Here, one operator and thus one
first integral was sufficient to work out the solution. In general,
one requires two, as we have a system of two first order equations,
which we wish to solve.

\subsection{One-Sector Model of Endogenous growth: The AK model}
We consider the following one-sector model of endogenous growth
presented in \cite{baro} where the representative consumer's utility
maximization problem is
 \be {Max} \quad \int_0^{\infty} e^{-(\rho-n) t}\frac{c^{1-\theta}-1}{1-\theta}dt ,\; \theta>0,\;\theta \not=1\label{(ck1)}\ee
subject to \be \dot a(t) =(r-n)a+w-c,\;c(0)=c_0\label{(cck2)},\ee
where $c(t)$ is the consumption per person, $a(t)$ is the assets per
person, $r(t)$ is the interest rate, $w(t)$ is the wage rate, and
$n$ is the growth rate of population. Suppose firms have the linear
production function \be y=f(k)=A k\ee where $A>0$. The marginal
product of capital is not diminishing, i.e. $f''=0$ and this
property makes it different from neoclassical production function.
The marginal product of capital is the constant $A$ and the marginal
product of labor is zero. Thus \be r=A-\delta,\ w=0 \ee where
$\delta \geq 0$ is the depreciation rate. It is assumed that the
economy is closed and $a(t)=k(t)$ holds. If we take $a=k$,
$r=A-\delta$ and $w=0$ then our optimal control problem is to
maximize {(\ref {(ck1)})} subject to \be \dot k =(A-\delta-n)k-c,
\;c(0)=c_0.\label{(ck2)}\ee

The current value Hamiltonian function is defined as \be
H(t,c,A,\lambda)=\frac{c^{1-\theta}-1}{1-\theta}+\lambda
[(A-\delta-n)k-c] \label{(ck3)},\ee where $c(t)$ is control
variable, $k(t)$ is the state variable and $\lambda(t)$  is the
costate variable. The necessary first order conditions for optimal
control are
 \be
\lambda=c^{-\theta} \label{(ck4)} \ee
 \be
\dot k =(A-\delta-n)k-c \label{(ck5)} \ee
 \be   \dot
\lambda+(A-\delta-\rho)\lambda =0.\label{(ck6)} \ee
 The
transversality condition is  \be \lim_{t\to\infty}
e^{-(\rho-n)t}\lambda(t)k(t)=0  \label{(aktrans)} \ee and from
{(\ref {(ck4)})} and {(\ref {(ck6)})}, the growth rate of
consumption is given by \be \frac{\dot
c}{c}=\frac{1}{\theta}(A-\delta-\rho) . \label{(akcon)} \ee

 Now we solve this model by utilizing our partial Hamiltonian approach. The partial Hamiltonian
operator determining equation  {(\ref{(12a)})}  for Hamiltonian
{(\ref{(ck3)})} with $\xi=\xi (t,k) $, $\eta=\eta (t,k) $,
$B=B(t,k)$ can be written as
\begin{eqnarray}
c^{-\theta}[\eta_t+((A-\delta-n)k-c)\eta_k]-\eta c^{-\theta}(A-\delta-n)\nonumber\\
-[\frac{c^{1-\theta}-1}{1-\theta}+c^{-\theta}((A-\delta-n)k-c)][\xi_t+((A-\delta-n)k-c)\xi_k]\nonumber\\
[-1.5ex]\label{(ck7)} \\[-1.5ex]
=B_t+B_k((A-\delta-n)k-c)-c^{-\theta}(\rho-n)
[\eta-\xi((A-\delta-n)k-c)],\nonumber
\end{eqnarray}
where we have used equations {(\ref {(ck4)})}-{(\ref
 {(ck6)})}. By following the same procedure for equation {(\ref{(ck7)})} as described in the previous
 two examples, we finally have
 \be \xi=a_1(t), \; \eta=[-\frac{\theta}{1-\theta}\dot a_1+(\rho-n)a_1]k+a_2(t),\;
  B=\frac{1}{1-\theta}a_1(t), \label{(ck8)}\ee
  \be -\frac{\theta}{1-\theta}\ddot a_1+[\rho-A+\delta-\frac{(\rho-n)\theta}{1-\theta}]\dot a_1
  +(\rho-n)(\rho-A+\delta)a_1=0,\label{(ck9)} \ee
\be \dot a_2+(\rho-A+\delta)a_2=0.\label{(ck10)}\ee Solving equation
{(\ref {(ck9)})} for $a_1(t)$, we have \be
a_1(t)=c_1e^{-(\rho-n)t}+c_2e^{\frac{(\rho-A+\delta)(1-\theta)}{\theta}t}
\label{(ck11)},\ee and equation {(\ref {(ck10)})} yields \be
a_2(t)=c_3e^{(A-\delta-\rho)t}. \label{(ck12)}\ee Thus $\xi, \eta,
B$ are given by
\begin{eqnarray}\xi=c_1e^{-(\rho-n)t}+c_2e^{\frac{(\rho-A+\delta)(1-\theta)}{\theta}t},
\nonumber\\
[-1.5ex]\label{(ckb13)} \\[-1.5ex]
 \eta=\frac{\rho-n }{1-\theta}c_1ke^{-(\rho-n)t}
-(n-A+\delta)c_2 ke^{\frac{(\rho-A+\delta)(1-\theta)}{\theta}t} +c_3e^{(A-\delta-\rho)t},\nonumber\\
  B=\frac{1}{1-\theta}[c_1e^{-(\rho-n)t}+c_2e^{\frac{(\rho-A+\delta)(1-\theta)}{\theta}t}], \nonumber \end{eqnarray}
 The following first integrals corresponding to operators
and gauge terms given in {(\ref {(ckb13)})} are computed from {(\ref
{(10)})}:
\begin{eqnarray}
I_1=e^{-(\rho-n)t}c^{-\theta} k [\frac{\rho-n
\theta+(\theta-1)(A-\delta)}{(1-\theta)}]
-\frac{\theta}{1-\theta}c^{1-\theta}e^{-(\rho-n)t}
\nonumber\\
[-1.5ex]\label{(ck13)} \\[-1.5ex] I_2=\frac{\theta}{\theta-1}c^{1-\theta}e^{\frac{(\rho-A+\delta)(1-\theta)}{\theta}t},\;
\; I_3=c^{-\theta} e^{-(\rho-A+\delta)t} . \nonumber\end{eqnarray}
Now we explain how to solve the AK model using the first integrals
$I_1, I_3$. Setting $I_1=a_1$ and $I_3=a_2$ after some
simplifications and using the initial condition $c(0)=c_0$, we have
\be k(t) =\frac{1-\theta}{\phi \theta}a_1
e^{(\rho-n)t}c^{\theta}(t)+\frac{1}{\phi} c(t) \label{(ck20)},\ee
where \be \phi=\frac{1}{\theta}[\rho-n
\theta+(\theta-1)(A-\delta)]\label{(ck21)}\ee and \be
c(t)=c_0e^{(\frac{A-\delta-\rho}{\theta})t},
c_0=a_2^{-\frac{1}{\theta}}.\label{(ck22)}\ee Note that
$(A-\delta-\rho)>0$ the consumption $c(t)$ given in {(\ref{(ck22)})}
increases with time. Substituting the value of consumption $c(t)$ in
{(\ref{(ck22)})} into {(\ref {(ck20)})}, capital is given by
 \be k(t)=\frac{1-\theta}{\phi \theta}a_1
c_0^{\theta} e^{(A-\delta-n)t}+\frac{1}{\phi}
c_0e^{(\frac{A-\delta-\rho}{\theta})t} ,\label{(ck23)}\ee where the
constant $a_1$ can be determined from the transversality condition.
The solutions {(\ref{(ck22)})} and {(\ref{(ck23)})} are the same as
given in \cite{baro} and here we deduced these by utilizing our
partial Hamiltonian approach. The transversality condition can be
rewritten as \be \lim_{t\to\infty} e^{-(\rho-n)t}c ^{-\theta}k=0
\label{(ck24)} \ee and we need to show $\lim_{t\to\infty}
e^{-(\rho-n)t}c ^{-\theta}k$ is zero. Using {(\ref{(ck22)})} and
{(\ref{(ck23)})} we have
  \be
\lim_{t\to\infty} \frac{1-\theta}{\phi \theta}a_1 +\frac{1}{\phi}
c_0^{1-\theta}e^{-\phi t} \label{(ck25)}\ee and it tends to zero
only if we choose constant $a_1=0$ and assume that
$\phi=\frac{1}{\theta}[\rho-n \theta+(\theta-1)(A-\delta)]>0$. This
further results in the following restriction on the parameters
 \be
\rho+\delta>(1-\theta)(A-\delta)+n \theta+\delta \label{(ck26)}\ee
and  {(\ref{(ck20)})} reduces to
 \be
k(t) =\frac{1}{\phi} c(t). \label{(ck27)}\ee

The detailed interpretations of the solutions obtained here are
given in \cite{baro}.

\section{Concluding remarks}
A systematic way to obtain reductions and closed-form solutions via
first integrals of Hamiltonian systems commonly arising in economic
growth theory and other economic models is developed.
 This is an algorithmic approach and can be applied to many state and
costate variables of the current value Hamiltonian. However, we
applied our method to systems with one control, one state and one
costate variable. The approach was explained with the help of one
simple illustrative example. We first studied two economic growth
models, the Ramsey model with a constant relative risk aversion
(CRRA) utility function and Cobb Douglas technology and the
one-sector AK model of endogenous growth. For the Ramsey model, the
solutions derived from our methodology were the same as those
derived by \cite{{baro},{runge}}. The restriction on the parameters
was obtained in a systematic way during the solution process unlike
in other models where it was assumed. The solutions were valid for
$\sigma \beta\not=1$, the case where the capital's share is not
equal to the intertemporal elasticity of substitution. The first
integrals and closed-form solutions for the one-sector AK model of
endogenous growth were also re-derived by our partial Hamiltonian
approach.

We  have shown that our systematic approach can be used to deduce
results given in the literature and we also found new solutions for
a variety of models.
\section*{Acknowledgments}
FMM is grateful to the staff of the Lahore School of Economics,
Pakistan for their warm hospitality during which time this work was
commenced and completed. In particular he thanks Dr. Shahid Amjad
Chaudhry for his friendship and constant support. FMM also thanks
the NRF of South Africa for research support through a grant.

\end{document}